\documentstyle[12pt,leqno,amstex]{article}

\input xypic

\begin{document}

\title{Hilbert-Kunz Functions of Cubic Curves and Surfaces}

\author{\small Ragnar-Olaf Buchweitz \\
        \small Department of Mathematics, University of Toronto\\ 
        \small 100 St. George Street \\ 
        \small Toronto, ON  Canada M5S 3G3 \\ 
        \small e-mail: ragnar@@math.toronto.edu 
        \and \small Qun Chen \\   
        \small Department of Mathematics, University of Toronto \\
        \small 100 St. George Street \\
        \small Toronto, ON  Canada M5S 3G3 \\
        \small e-mail: chen@@math.toronto.edu} 
\maketitle

\begin{abstract}

We determine the Hilbert-Kunz function of plane 
elliptic curves in odd characteristic, as well as 
over arbitrary fields the generalized Hilbert-Kunz 
functions of nodal cubic curves.

Together with results of K.~Pardue and P.~Monsky, this
completes the list of Hilbert-Kunz functions
of plane cubics. 

Combining these results with the calculation of the
(generalized) Hilbert-Kunz function  of Cayley's
cubic surface, it follows that in each degree 
and over any field of positive characteristic
there are curves resp. surfaces taking on the minimally 
possible Hilbert-Kunz multiplicity.

\end{abstract}

\vfill



\section{Introduction}
Let $S=k[x_0,\cdots,x_n]$ be the standard polynomial ring in
$n+1$ variables over a field $k$ of prime characteristic
$p$. Given a finite graded $S$-module $M$, the Hilbert-Kunz
function of $M$ is defined on powers of the characteristic,
$q=p^n$, $n\in {\Bbb{N}}$, through 
\[ HK_M (q):=\dim_k M/m^{[q]}M \]
where $m^{[q]}=(x^q_0,\cdots,x^q_n)$ is the
$q$-th Frobenius power of the maximal homogeneous ideal 
$m=(x_0, \cdots, x_n)$.
If $I\subset S$ is a homogeneous ideal, and
$R=S/I$ the homogeneous coordinate ring of the 
underlying projective scheme $X\subset {\Bbb{P}}^n_k$, 
the function $HK_R(q)$ is also called the
Hilbert-Kunz function of $X$.
Introduced by E. Kunz [5] in 1969, 
these functions were first studied in detail by 
P. Monsky [6], and he obtained the
asymptotic formula
\[ HK_M (q)=c q^m+ O(q^{m-1})  \]
with $c\geq 1$ some real number and 
$m$ the Krull dimension of $M$. 
The number $c$ is called the Hilbert-Kunz multiplicity of
$M$ and P. Monsky conjectures it to be rational. 
In general, it seems very difficult to determine these
functions explicitly and a conceptual interpretation of the 
constant $c$ is missing 
(see [4] for some surprising examples). 

Here we exhibit the Hilbert-Kunz functions of plane elliptic
curves in odd characteristic and of plane nodal cubics.
Combining this work with results in [7] and [8]
completes the explicit determination of Hilbert-Kunz
functions of plane cubic curves. The Hilbert-Kunz functions
of reducible cubics were already determined by K. Pardue in
[8], and he also predicted the following list for the
irreducible ones on the basis of computer experiments. Note
that the Hilbert-Kunz function is invariant under extensions
of the coefficient field $k$, so that one may assume $k$
algebraically closed.

\newtheorem{theorem}{Theorem}
\begin{theorem}
Let $f$ be the equation of an irreducible cubic curve 
$C$ in ${\Bbb{P}}^2_k$ over an algebraically closed field 
$k$, and let $HK_R(q)$ be the Hilbert-Kunz function of 
the homogeneous coordinate ring $R=S/f$.
\begin{description}
\item[(1) (K. Pardue {[8]})] If $C$ is a cuspidal cubic,
\[ HK_R (q)=\left\{            
\begin{array}{ll}            
\frac{7}{3} q^2 &\text{for}\quad p=3\; , \\            
\frac{7}{3} q^2-\frac{4}{3} & \text{for}\quad p\neq 3\; .           
\end{array}             
\right. \]
\item[(2) (Thm. \ref{th:nodal} below)] If $C$ is a nodal cubic,
\[ HK_R (q)=\left\{ 
\begin{array}{ll}     
\frac{7}{3}q^2-\frac{1}{3}q-1 &      
	\text{for}\quad q\not\equiv 2\bmod 3\;,  \\     
\frac{7}{3}q^2-\frac{1}{3}q-\frac{5}{3} &     
	\text{for}\quad q\equiv 2\bmod 3\;.     
\end{array}     
\right. \]
\item[(3) (Thm. \ref{theorem:main} below)] 
If $C$ is an elliptic curve and $p\neq 2$,
\[ HK_R(q)=\frac{9}{4} q^2-\frac{5}{4}\quad .\]
\item[(4) (P. Monsky {[7]})] If $C$ is an elliptic curve 
and $p=2$,
\[ HK_R (q)=\left\{    
\begin{array}{ll}   
\frac{9}{4}q^2 & \text{if the $j$-invariant is $0$}\;,  \\
\frac{9}{4}q^2-1 &\text{if the $j$-invariant is not $0$}\;.    
\end{array}   
\right. \]
\end{description}
\label{theorem:major}
\end{theorem}

At this stage, Hilbert-Kunz functions or multiplicities of 
plane curves of higher degree remain mysterious. 
However, a corollary of our work shows that for any 
$d\ge 2$ and for any field $k$ of 
prime characteristic there exists a plane curve of 
degree $d$ in ${\Bbb{P}}^2_k$ whose 
Hilbert-Kunz multiplicity is $\frac{3}{4} d$ 
--- and this is the minimal possible value
for such curves.
In particular, the minimal Hilbert-Kunz multiplicity in 
each degree is rational and 
independent of the characteristic. 
We then determine explicitly the Hilbert-Kunz function 
of Cayley's cubic surface in ${\Bbb{P}}^3_k$, 
and the result allows us to conclude as well 
that for any $d\ge 2$ and for any field $k$ 
of prime characteristic there exists a surface 
of degree $d$ in ${\Bbb{P}}^3_k$ whose  
Hilbert-Kunz multiplicity is $\frac{2}{3} d$ --- 
and this is again the minimal possible value, 
again rational and independent
of the characteristic. 

\section{Minimal Values of Hilbert-Kunz Functions}
Let $I=(f)$ be a principal ideal generated by a 
homogeneous form $f$ of degree $d > 0$ in $S$. 
The considerations in this section apply to the 
values of the {\em generalized\/} Hilbert-Kunz 
function of $R=S/I$, introduced by A. Conca in [1], and
defined as
\[ HK_{R,x}(q) =\dim_k S/(f,x_0^{q},\ldots,x_n^{q}) \quad, \]
where $q$ is now {\em any\/} nonnegative integer, 
$k$ any field. 
Unless $k$ is of positive characteristic $p$, 
and $q$ is a power of $p$,
this dimension will generally depend upon the choice 
of the coordinate system 
$x = (x_0,\ldots,x_n)$.

For each $q\in{\Bbb N}$ and each choice of coordinates $x$,
set $x^{[q]}=(x_0^{q},\ldots,x_n^{q})$ and
consider the following graded $S$-modules of finite length,
\begin{eqnarray*}
\Theta & = & \frac{S}{x^{[q]}}=\bigoplus_i \Theta_i  \quad , \\ 
\theta & = & \frac{S}{f+x^{[q]}}=\bigoplus_i \theta_i \quad ,  \\
\vartheta & = & \frac{(x^{[q]}:f)}{x^{[q]}}=\bigoplus_i \vartheta_i \quad . 
\end{eqnarray*}
They are related by the exact sequence of graded 
$S$-modules 
\begin{equation} 
\begin{array}{cccccccccccc}  0 &
\longrightarrow & \vartheta (-d)&  \longrightarrow & \Theta (-d) &  
\stackrel{f}{\longrightarrow} & \Theta &   \longrightarrow & \theta &  
\longrightarrow & 0 &,   
\end{array}
\label{eq:exa}
\end{equation}
and $HK_{R,x}(q)=\dim_k \theta$. 
Evaluating dimensions yields universal bounds for the 
generalized Hilbert-Kunz function of 
$R=S/(f)$, when $f$ varies over polynomials of degree 
$d$ in $n+1$ variables,
\begin{equation}
q^{n+1}=\sum_i \dim_k \Theta_i \geq HK_{R,x}(q)\geq
\sum_i \max    \{ \dim_k \Theta_i-\dim_k \Theta_{i-d}, 0 \} \quad.
\label{eq:ineq}
\end{equation}
The upper bound, $HK_{R,x}(q)=q^{n+1}$, is achieved iff 
$f\in x^{[q]}$; 
for example if $d > (n+1)(q-1)$, 
or if $q=p^n$ is a power of the characteristic, 
$d\geq q$ and $f=l^d$ for some linear form $l$.
Here we are more concerned with the lower bound, 
that is taken on
if and only if $f$ {\em is of maximal rank at\/} $q$, 
meaning that in each degree $i$ the $k$-linear map 
$f|\Theta_{i-d}$ is of maximal rank. 

Whether a given polynomial $f$ is of maximal rank at $q$,
can be decided by looking at the socle degree of the 
artinian ring $\theta$,
\begin{equation}a(q) =\max\{i:\theta_i\neq 0\}\quad, 
\label{eq:a}
\end{equation}
and at the initial degree of $\vartheta$, 
\begin{equation}
\imath(q)=\min\{i:\vartheta_i \neq 0\}\quad.
\label{eq:i}
\end{equation}
Indeed, as the socle degree of $\Theta$ is $(n+1)(q-1)$, 
outside the range $d\leq i\leq (n+1)(q-1)$ source or 
target of $f|\Theta_{i-d}$ is zero, whereas for
a degree $i$ inside that range the map is not surjective 
iff $i \leq a(q)$, not injective 
iff $i-d \geq \imath(q)$. 
Accordingly, all the $k$-linear maps induced by $f$ 
are of maximal rank 
iff $a(q) < \imath(q) + d$.
Moreover, the exact sequence (\ref{eq:exa}) is selfdual, 
whence it suffices
to know either $a(q)$ or $\imath(q)$:

\newtheorem{lemma}{Lemma}
\begin{lemma}
For each $q\in{\Bbb N}$, and independent of $f$, one has
\begin{equation}
a(q)+\imath(q)=(n+1)(q-1)\quad .
\end{equation}
Given $q$, all $k$-linear maps $f|\Theta_{i-d}$ are 
of maximal rank iff
\begin{equation}
a(q) < \frac{(n+1)(q-1)+d}{2} < \imath(q)+d\quad .
\end{equation}
Moreover, each of the inequalities implies the other.
\label{lemma:mr}
\end{lemma}
{\noindent \it Proof:\quad }  
The ring $\Theta=S/x^{[q]}$ is a zerodimensional 
complete intersection with its
socle in degree $(n+1)(q-1)$.
Thus for any finite graded $\Theta$-module $M$, 
we have an
isomorphism of graded $\Theta$-modules
\[ {\rm Hom}_k (M,k) \cong {\rm Hom}_{\Theta}
(M,\omega_{\Theta})\quad ,   \]
where $\omega_{\Theta}=\Theta ((n+1)(q-1))$ is 
the canonical module of $\Theta$,  
and ${\rm Hom}_k (M,k)$ is the $\Theta$-module graded 
naturally through
\[ ({\rm Hom}_k (M,k))_i = {\rm Hom}_k (M_{-i},k)\quad .  \] 
As ${\rm Hom}_{\Theta}(\theta,\Theta)   
\cong (x^{[q]}:f)/x^{[q]}=\vartheta$, 
we get
\[ {\rm Hom}_k (\theta,k)\cong \vartheta ((n+1)(q-1))\quad . \]
For the dimension of the finite dimensional 
$k$-vector space $\theta_i$, this yields  
\[ \dim_k \theta_i =\dim_k {\rm Hom}_k (\theta_i,k)
       = \dim_k ({\rm Hom}_k (\theta,k))_{-i}
       = \dim_k \vartheta_{(n+1)(q-1)-i}\quad,  \]      
and the equality follows from the definition of 
$a(q)$ and $\imath(q)$.
As $f$ induces maps of maximal rank iff $a(q) < \imath(q) + d$, 
we can eliminate either one of the two invariants to obtain 
the last claim.
$\hfill \square$ 
\[ \]

If $d>(n+1)(q-1)$, the information is already complete:
$f$, inducing the zero map in (\ref{eq:exa}),
is trivially of maximal rank at $q$, and 
$HK_{R,x}(q) = q^{n+1}$.
Also, if $f$ is a polynomial of a single variable, $n=0$, 
there are no secrets to discover. 
If $n>0$, the (usual) Hilbert series 
\[  H_{\Theta}(t) = \sum_i (\dim_k \Theta_i)t^i = 
(1+t+t^2+\cdots+t^{q-1})^{n+1}\quad,\] 
of the artinian $k$-algebra 
$\Theta$ is a {\em reciprocal\/} and {\em unimodal\/} 
polynomial of degree $l=(n+1)(q-1)$ in $t$, meaning that
its coefficients, $\alpha_i=\dim_k \Theta_i$, satisfy
\begin{eqnarray*}
\alpha_i & = & \alpha_{l-i} \quad {\rm for\ every}\quad i \quad, \\
\alpha_i & < & \alpha_{i+1} \quad {\rm for} \quad 
 0\leq i <\left\lfloor \frac{l}{2} \right\rfloor\quad .
\end{eqnarray*}
In particular, 
$\dim_k\Theta_i - \dim_k\Theta_{i-d} > 0$ iff
$0\leq i \leq m(q)$, where
\[ m(q) =\left\lfloor \frac{(n+1)(q-1)+(d-1)}{2} 
		 	\right\rfloor\quad,\]
--- as for $a(q)$, we suppress the dependence upon 
$d$ from the notation.
Thus the lower bound, $L(q)$, in inequality (\ref{eq:ineq}) 
evaluates to
\begin{eqnarray*}
   L(q) & := & \sum_i \max\{ \dim_k \Theta_i-\dim_k 
                        \Theta_{i-d}, 0 \} \\
    & = & \sum_{m(q)-d+1}^{m(q)} \dim_k \Theta_i 
    \quad,\quad {\rm as}\ H_{\Theta}(t)\ {\rm is\ unimodal\ 
                           and\ reciprocal}\\
    & = & {\rm\ coefficient\ of\ \;} t^{m(q)} {\rm\;\ in\ }
          \frac{(1-t^d)(1-t^q)^{n+1}}{(1-t)^{n+2}}\\
    & = & \frac{1}{2\pi\sqrt{-1}}\int_{|z|=\epsilon}
\frac{(1-z^d)(1-z^q)^{n+1}}{(1-z)^{n+2}z^{m(q)+1}}dz\quad.
\end{eqnarray*}

\noindent
As $m(q)$ is the largest integer smaller than 
$\frac{1}{2}((n+1)(q-1)+d)$,
we get the following result.

\begin{theorem}  
If $n>0$, and if $f$ is a homogeneous polynomial 
of degree $d\leq (n+1)(q-1)$\; in \;$n+1$\; many variables, 
then the socle degree of the graded artinian $k$-algebra 
$\theta$ satisfies
\begin{equation}
a(q) \geq m(q)\quad .
\label{eq:aq}
\end{equation}
Furthermore, the following statements are equivalent:
\begin{description}
\item[{\rm (i)}]\quad 
	The polynomial $f$ is of maximal rank at $q$.
\item[{\rm (ii)}]\quad 
	The socle degree $a(q)$ is minimal, $a(q)=m(q)$.
\item[{\rm (iii)}]\quad 
	The initial degree $\imath(q)$ is maximal, 
	$\imath(q)=(n+1)(q-1)-m(q)$.
\item[{\rm (iv)}]\quad 
	The Hilbert-Kunz function of $f$ at $q$ achieves 
	the lower bound $L(q)$.
\end{description}
\label{theorem:minimal}
\end{theorem}

{\noindent \it Proof:\quad}
The first statement follows from the exact sequence 
(\ref{eq:exa}), as 
\[ \dim_k \Theta_{m(q)} > \dim_k \Theta_{m(q)-d}\quad . \]

As $a(q)$ is an integer and $m(q)$ is the largest integer 
smaller than $\frac{1}{2}((n+1)(q-1)+d)$, the just 
established lower bound for $a(q)$ implies the equivalences 
in view of Lemma \ref{lemma:mr}. 
$\hfill \square$ 
\[ \]
 
{\noindent \it Example:\quad} For $d=2,3$ and $n=2,3$, 
we get the following table:
\[ \begin{tabular}{|l|l|l|l|} \hline 
   \emph{$d$} & \emph{$n$} & \emph{$m(q)$} & 
   \emph{lower bound for $HK_{R,x} (q)$}  \\ \hline
  2 & 2 & $\left\lfloor \frac{3q}{2} \right\rfloor -1$ & 
  $\left\lfloor \frac{3}{2} q^2 \right\rfloor 
  =\left\{ \begin{array}{ll}
         \frac{3}{2} q^2  & \mbox{for $q$ even} \\ 
         \frac{3}{2} q^2-\frac{1}{2} & \mbox{for $q$ odd}
         \end{array}
\right.$  \\ \cline{2-4} 
    & 3 & $2(q-1)$ & $\frac{4}{3} q^3-\frac{1}{3} q$  \\ \hline
  3 & 2 & $\left\lfloor \frac{3q-1}{2} \right\rfloor$ &
  $\left\lfloor \frac{9}{4} q^2-\frac{5}{4} \right\rfloor
  =\left\{ \begin{array}{ll}
         \frac{9}{4} q^2-2 & \mbox{for $q$ even} \\  
         \frac{9}{4} q^2-\frac{5}{4} & \mbox{for $q$ odd}
         \end{array}
\right.$  \\ \cline{2-4} 
    & 3 & $2q-1$ & $2q^3-q$  \\ \hline 
  \end{tabular}
\]   
For $d=2$, it can be extracted from [1] that 
the quadric $x_0^2-x_1x_2$ for $n=2$, respectively
the quadric $x_0x_1-x_2x_3$ for $n=3$,
have generalized Hilbert-Kunz functions that 
take on the minimum value at each $q$. 
\[ \]

{\noindent \it Remark 1:\quad} 
P. Monsky noted that expressing the minimal 
possible value $L(q)$ of $HK_{R,x}(q)$ as a residue 
leads to an intriguing lower bound for 
Hilbert-Kunz multiplicities in terms of the
integrals 
\[ \beta_{n+1}=\frac{1}{\pi} \int^{+\infty}_{-\infty}
(\frac{\sin \alpha}{\alpha})^{n+1} d\alpha = 
\frac{1}{2^{n}n!}\sum_{i=0}^{\left\lfloor n/2 \right\rfloor}
(-1)^i(n+1-2i)^{n}\binom{n+1}{i}\quad , \]
as
\begin{equation}
d\beta_{n+1} = \lim_{q\to\infty}\frac{1}{q^n{2\pi\sqrt{-1}}}
\int_{|z|=\epsilon}
\frac{(1-z^d)(1-z^q)^{n+1}}{(1-z)^{n+2}z^{m(q)+1}}dz\quad.
\label{eq:beta}
\end{equation}
Thus, for a hypersurface of degree $d$ in ${\Bbb{P}}^n_k$ 
over a field $k$ of positive characteristic, 
the Hilbert-Kunz multiplicity satisfies
\begin{equation}
c \geq d\beta_{n+1}\quad.
\label{eq:lower}
\end{equation}
A direct combinatorial proof is as follows. Expanding 
$(1-t^d)/(1-t)^{n+2}$ into its Taylor series at $t=0$, 
one can write
\[
\frac{1-t^d}{(1-t)^{n+2}} = 
\sum_{\nu=0}^{d-n-1} R(\nu)t^{\nu} +
\sum_{\nu \geq 0} P(\nu)t^{\nu}\quad,\]
where $R(\nu)\in{\Bbb Z}$, and 
$P(\tau) = \frac{d}{n!}\tau^n + O(\tau^{n-1})$ is the 
corresponding (Hilbert) polynomial, univariate over ${\Bbb Q}$
of degree $n$ with leading coefficient $\frac{d}{n!}$.
Now use that $m(q)/q = (n+1)/2 + O(1/q)$, that the
coefficient of $t^{m(q)}$ in 
$\frac{1}{q^n}(1-t^q)^{n+1}\sum_{\nu=0}^{d-n-1} R(\nu)t^{\nu}$
tends to zero with $q$, and that
\[
(1-t^q)^{n+1}\sum_{\nu \geq 0} P(\nu)t^{\nu} =
\sum_{\nu\geq 0}
\big(\sum_{i=0}^{\left\lfloor \nu/q \right\rfloor}
(-1)^i P(\nu - iq)\binom{n+1}{i}\big)t^{\nu}\quad,
\]
to get
\begin{eqnarray*}
\lim_{q\to\infty} \frac{L(q)}{q^n} &=&
\lim_{q\to\infty} 
	\sum_{i=0}^{\left\lfloor m(q)/q\right\rfloor}
		(-1)^i\frac{P(m(q)-iq)}{q^n}\binom{n+1}{i}\\
&=&\lim_{q\to\infty} 
	\sum_{i=0}^{\left\lfloor m(q)/q\right\rfloor}
		(-1)^i\frac{d}{n!}(\frac{m(q)}{q}-i)^n\binom{n+1}{i}\\
&=&\frac{d}{2^n n!}
\sum_{i=0}^{\left\lfloor n/2\right\rfloor}
(-1)^i(n+1-2i)^n\binom{n+1}{i}\\
&=&d\beta_{n+1}\quad.
\end{eqnarray*}
That this combinatorial expression equals the indicated integral
can now be checked in any table of integrals, e.g. 
[3, $3.836.5^3$, p 458]. It follows that the sequence
$\{ \beta_{n} \}$ of rational numbers decreases to zero. 
The first few values are
\[ \beta_1=1, \; \;\beta_2=1, \; \;  
\beta_3=\frac{3}{4}, \; \;      \beta_4=\frac{2}{3}, \; \;    
\beta_5=\frac{115}{192}, \; \;
\beta_6=\frac{11}{20}\quad .  \]

{\noindent \it Remark 2:\quad} 
The same argument applies to $\liminf_{q\to\infty}
HK_{R,x}(q)/q^n$, but it is not yet known whether a 
``generalized'' Hilbert-Kunz multiplicity exists, i.e.
whether
$HK_{R,x}(q)/q^n$ tends to a limit for $q\to\infty$.
\[ \]

For elliptic curves in odd characteristic, we will prove 
that the corresponding cubic polynomial is of maximal rank 
at any power $q$ of the characteristic, 
whereas for the polynomial $x_0x_1x_2x_3(
\frac{1}{x_0}+\frac{1}{x_1}+\frac{1}{x_2}+\frac{1}{x_3})$, 
representing Cayley's cubic surface, this will be even 
established at any $q\in{\Bbb N}$
over any field $k$.
The proof is accomplished by showing that $a(q)$ equals the 
minimum value $m(q)$, and the respective Hilbert-Kunz function 
can then be read off from the table above.

None of this applies to elliptic curves in characteristic $2$, 
(see [7] for details), nor does it hold for singular 
irreducible cubic curves in any characteristic. 
But in the latter case, the Hilbert-Kunz function can be 
determined completely from the rational parametrization 
of the curve as we show next.

\section{Singular Irreducible Cubic Curves}
The Hilbert-Kunz function of a cuspidal cubic is known from [8], 
see also [1], but our treatment here deals with the nodal 
and cuspidal case at the same time. 
Let $k$ be an algebraically closed field 
--- of any characteristic for now ---, 
and denote by $C$ a singular irreducible 
plane cubic curve in ${\Bbb{P}}^2_k$. 

In suitable coordinates, 
$C$ is given by a Weierstra\ss\ equation
$$
f(x,y,z) = z(y^2 + a_1 xy - a_2 x^2) - x^3 = 0 \quad,
$$
so that $o=[0,0,1]\in{\Bbb{P}}^2$ is its unique singular point. 
The curve has a node at $o$ iff the tangential quadric
$Q(x,y) = y^2 + a_1 xy - a_2 x^2$ has distinct roots iff
$a^2_1 + 4 a_2 \ne 0$, otherwise it is cuspidal.

The curve $C$ is rational and a rational parametrization
$\nu:{\Bbb{P}}^1\to C\subset {\Bbb{P}}^2$ normalizes the curve,
pulling back ${\cal{O}}_{{\Bbb{P}}^2}(1)$ along 
$C\hookrightarrow{\Bbb{P}}^2$ and then $\nu$ to ${\cal{O}}_{{\Bbb{P}}^1}(3)$. 
Algebraically, such a parametrization is given by the monomorphism of
$k$-algebras 
\[ \alpha(x,y,z)=(sQ(s,t),tQ(s,t),s^3)\quad , \]
$$
\alpha: R ={ k[x,y,z]\over f(x,y,z)} \cong 
\bigoplus_{n\in {\Bbb{Z}}}
H^0(C,{\cal{O}}_{C}(n))\hookrightarrow 
\bigoplus_{n\in {\Bbb{Z}}}
H^0({\Bbb{P}}^1, {\cal{O}}_{{\Bbb{P}}^1}(3n)) 
\cong k[s,t]^{(3)} =:{\tilde R}\quad,
$$
where ${\tilde R}=k[s,t]^{(3)}$ is the Veronese subring 
of the polynomial ring $k[s,t]$ spanned by all homogeneous 
polynomials whose degree is divisible by $3$. 
Notice that ${\tilde R}_n$ consists of all homogeneous 
polynomials of degree $3n$.

The cokernel of $\alpha$ can be identified as follows. 
A section $p(s,t)\in H^0({\Bbb{P}}^1, {\cal{O}}_{{\Bbb{P}}^1}(3n))$ 
comes via $\alpha$ from a section in
$H^0(C,{\cal{O}}_{C}(n))$ iff $p(s,t)$ takes on the 
same value at the two points $Q(s,t)=0$. 
Explicitly, write $p(s,t) = e_1(s) + e_2(s) t + e(s,t)Q(s,t)$
with uniquely determined polynomials $e_1,e_2\in k[s]$ 
and $e\in k[s,t]$.
The component $e_2(s)t$ represents the class of 
$p(s,t)$ in $k[s,t]/(k[s] + Q k[s,t])$ and 
$(k[s,t]/(k[s] + Q k[s,t]))^{(3)}$ is the cokernel of $\alpha$.
If $p(s,t)\in k[s,t]^{(3)}$, then $e_2(s) = b(s^3) s^2$ 
for some unique univariate polynomial $b$ that is necessarily 
of degree $\frac{1}{3}\deg p -1$.

\begin{lemma}
{\rm (i)}\quad
The map
$$
\beta : {\tilde R} = k[s,t]^{(3)}\to k[z](-1)\quad,
$$
associating to $p(s,t)$ the polynomial $b(z)$, 
is a degree preserving epimorphism of $R$-modules, 
the $R$-module structure on ${\tilde R}$ given by
$\alpha$, the one on $k[z]$ by the natural projection 
$R=k[x,y,z]/f\to k[z]$.

\noindent {\rm (ii)}\quad
The sequence of graded $R$-modules
\begin{equation}
\begin{array}{ccccccccc}
 0 & \to & R & \stackrel{\alpha}{\longrightarrow} &
  {\tilde R} = k[s,t]^{(3)} & 
  \stackrel{\beta}{\longrightarrow} & k[z](-1) &
   \to & 0
   \end{array}
   \label{eq:exa2}
   \end{equation}
   is exact.
   \label{lemma:normalization}
   \end{lemma}

{\noindent \it Proof:\quad }
If $p(s,t) = e_1(s) + b(s^3)s^2 t + e(s,t)Q(s,t)$ is the unique
representation of $p(s,t)\in {\tilde R}$, then
\begin{eqnarray*} 
\alpha(x)p(s,t) &=& 0 + 0\cdot s^2t + (sp(s,t))Q(s,t) \quad, \\
\alpha(y)p(s,t) &=& 0 + 0\cdot s^2t + (tp(s,t))Q(s,t) \quad, \\
\alpha(z)p(s,t) &=& s^3e_1(s) + (s^3b(s^3))s^2 t + (s^3e(s,t))Q(s,t) \quad,
\end{eqnarray*}
are the corresponding unique representations of 
$\alpha(x)p(s,t), \alpha(y)p(s,t)$ and $\alpha(z)p(s,t)$ 
respectively. 
This shows that the image of $\beta$ is annihilated by 
$x,y$ and that $\beta(\alpha(z)p) = z\beta(p)$. 
Furthermore, $\beta(s^2t)$ generates
the image of $\beta$ already as $k[z]$-module, 
thus a fortiori as $R$-module. 
As $s^2t$ is of degree one with respect to the grading
on ${\tilde R} = k[s,t]^{(3)}$, (i) follows.

For (ii), note first that $\beta\alpha(1) = 0$,
whence $\beta\alpha = 0$. 
To prove that the kernel of $\beta$ is precisely the
image of $\alpha$, consider Hilbert functions: 
In degree $i\in {\Bbb N}$,
\[ \dim_k R_i =\left\{ 
   \begin{array}{ll}
   1  & for \quad i=0 \\ 
   3i & for \quad i>0 
   \end{array}
   \right.
   \quad,
   \quad {\rm whereas}
   \quad \dim_k {\tilde R}_i = 3i + 1\quad. \]
Accordingly, the quotient ${\tilde R}_i/R_i$ is zero 
for $i=0$ and onedimensional for $i>0$. 
Thus the cokernel of $\alpha$ and $k[z](-1)$ have
the same Hilbert function and (ii) follows.
$\hfill\square$
\[ \]

Multiplication with $x^q,y^q,z^q$ on (\ref{eq:exa2}) 
results in a commutative diagram of graded $R$-modules 
whose exact rows and columns define the modules
$A$ through $G$,
\[ 
\diagram
 & 0 \dto & 0 \dto & 0 \dto  \\
0 \rto & A  \dto \rto^{\matrix \overline{\alpha} \endmatrix} & 
B \dto \rto^{\matrix \overline{\beta} \endmatrix} &
C \dto \rto & D \rto & 0 \\
0 \rto & R^{\oplus 3}(-q) \xto[ddd]_{\left(
                                 \matrix 
                                 x^q  \\
                                 y^q \\
                                 z^q
                                 \endmatrix
                                 \right)} 
                          \rto^-{\matrix \alpha^{\oplus 3} \endmatrix} &
{\tilde R}^{\oplus 3}(-q) \xto[ddd]_{\left(
                                 \matrix
                                      \alpha(x)^q \\
                                      \alpha(y)^q \\
                                      \alpha(z)^q
                                 \endmatrix
                                 \right) }
                           \rto^-{\matrix \beta^{\oplus 3} \endmatrix} &
k[z]^{\oplus 3}(-1-q) \xto[ddd]_{\left(
                             \matrix
                                   0 \\
                                   0 \\
                                   z^q
                             \endmatrix
                             \right)} \rto & 0 \\
&&&&                          \\
&&&&                          \\
0 \rto & R \dto \rto^{\matrix \alpha \endmatrix} & 
{\tilde R} \dto \rto^{\matrix \beta \endmatrix} &
k[z](-1) \dto \rto & 0 \\ 
&  E \dto \rto & F \dto \rto & G \dto \rto & 0 \\
&  0 & 0 & 0 &&\quad. 
\enddiagram 
\]
In this diagram, $E,F,G$, and then also $D$, 
are finite dimensional and one has
\begin{equation} 
HK_{R,(x,y,z)}(q) = \dim_kE = \dim_k F - \dim_k G + \dim_k D\quad.
\label{eq:sum}
\end{equation}

The dimension of $G\cong (k[z]/z^q)(-1)$ equals $q$, and 
the next Lemma determines the dimension of $F$, that is the value
of the generalized Hilbert-Kunz function for ${\tilde R}$
with respect to $(x,y,z)$ at $q$.

\begin{lemma} 
{\rm (i)}\quad
Set $P=k[s,t]$, the polynomial ring in two variables 
with its natural grading.
For any $q\in{\Bbb{N}}$, the $P$-module 
$M = k[s,t]/(\alpha(x)^q,\alpha(y)^q,\alpha(z)^q)$ 
has minimal graded resolution
$$ \diagram 
   0 \rto & 
   P(-4q)\oplus P(-5q)  
   \xto[rrrr]^-{\left(
    \matrix
     t^q & -s^q & 0 \\
  s^{2q} & 0    &-Q(s,t)^q
    \endmatrix
    \right) } & & & & 
  P(-3q)^{\oplus3}  
    \xto[rr]^-{\left(
         \matrix 
         s^qQ(s,t)^q \\
         t^qQ(s,t)^q \\
         s^{3q}
         \endmatrix
         \right)}   & &   
  P \rto & M \rto   &  0\ . 
\enddiagram
$$
{\rm (ii)}\quad
The middle column in the diagram above is obtained from that
resolution by applying the functor $(\ \;)^{(3)}$, and in 
particular $F=M^{(3)}$.\\
\noindent {\rm (iii)}\quad
\[ HK_{{\tilde R},(x,y,z)}(q) = \dim_k F =
        \left\{
        \begin{array}{ll} 
	{7\over 3} q^2 - {1\over 3} & \text{if}\quad q\not\equiv 0\bmod 3\quad, \\
	{7\over 3} q^2 & \text{if} \quad q\equiv 0\bmod 3\quad.
        \end{array}
        \right. 
\]
\label{lemma:F}
\end{lemma}

{\noindent \it Proof:\quad } 
(i)\quad As $s$ does not divide $Q(s,t)$, the module $M$ is artinian,
and the result follows from the Hilbert-Burch theorem: 
the (signed) $(2\times 2)$-minors of the
leftmost matrix are respectively
$s^q Q(s,t)^q = \alpha(x)^q\; ,\; t^q Q(s,t)^q = \alpha(y)^q\;,
\; s^{3q} =\alpha(z)^q\;$.

\noindent (ii) is clear and (iii) follows then easily from
\[ \dim_k F_i = \dim_k M_{3 i} = 
\dim_k P_{3 i} - 3 \dim_k P_{3 i - 3 q} + \dim_k P_{3 i - 4 q} +
\dim_k P_{3 i - 5 q}\]
and $\dim_k P_j = \max\{0,j+1\}$ for $j\in\Bbb Z$.
$\hfill\square$
\[ \]

In equation (\ref{eq:sum}) for $HK_{R,(x,y,z)}(q)$, 
it remains to determine the dimension of $D$.
To this end, we exhibit the map ${\overline \beta}$ explicitly. 
As multiplication by $z^q$ is injective on $k[z]$,
one has $C\cong k[z]^{\oplus 2}(-1-q)$. 
Furthermore, as $B \cong \big( k[s,t](-4q)\oplus k[s,t](-5q)\big)^{(3)}$ 
by Lemma \ref{lemma:F}, a homogeneous element in 
$B$ is represented by a pair $(p_1(s,t),p_2(s,t))$ 
of homogeneous polynomials satisfying 
\[ \deg p_1 = \deg p_2 + q \equiv -4q\bmod 3\quad,\]
and such pair is mapped to
\[
p_1(s,t)(t^q,-s^q,0) + p_2(s,t)(s^{2q},0,-Q(s,t)^q)
\]
in ${\tilde R}^{\oplus 3}(-q)$. Thus 
\[
{\overline \beta}(p_1,p_2) = (\beta(t^qp_1 + s^{2q}p_2),\beta(-s^q p_1))
\in k[z]^{\oplus 2}(-1-q)\cong C\quad.
\]

As the field $k$ is algebraically closed, the quadric $Q$ 
factors,
$Q(s,t) = (t-us)(t-vs)\;;\; u,v\in k$;
and the unique representation of $t^q \bmod Q(s,t)$ is
\[
t^q = \tau_1 s^q + \tau_2 s^{q-1}t + \tau(s,t) Q(s,t)\quad,
\]
where
\[
\tau_1 =   -uv\sum_{i=0}^{q-2} u^{q-2-i}v^{i} \quad,\quad
	\tau_2 = \sum_{i=0}^{q-1} u^{q-1-i}v^{i}\quad.
\]
Writing now 
\[
p_1 = b_1(s) + b_2(s)t + b_3(s,t)Q(s,t)\quad,\quad
	p_2 = c_1(s) + c_2(s)t + c_3(s,t)Q(s,t)\quad,
	\]
for suitable polynomials $b_i,c_i$, it follows that
\[
  {\overline \beta} (p_1,p_2) = 
\big(\tau_2 b_1 s^{q-3}+ (\tau_1 - a_1 \tau_2) b_2 s^{q-2} + 
c_2 s^{2q-2}, -b_2 s^{q-2}\big)\big|_{s^3=z}\quad.
\]
So the image of ${\overline \beta}$ in 
$C\cong k[z]^{\oplus 2}(-1-q)$ is generated as a $k[z]$-module
by the three pairs
\[ (\tau_2,0)z^{(q-3+\epsilon)/3}\quad,\quad
(\tau_1 -a_1 \tau_2 ,-1)z^{(q-2+\eta)/3}\quad,\quad
(1,0)z^{(2q-2+\zeta)/3}\quad,
\]
where $\epsilon,\eta,\zeta\in \{0,1,2\}$ are such that
the exponents become integers. Accordingly, the dimension of $D$
is equal to $(2q-5 +\epsilon+\eta)/3$ if $\tau_2\ne 0$, 
whereas it equals $(3q-4+\eta+\zeta)/3$ if $\tau_2 = 0$. 
Thus, the decisive
factor is whether or not $\tau_2$ vanishes.

\begin{lemma} 
\begin{description}
\item[(i)]
If $q$ is a power of the characteristic of $k$, 
then $C$ is nodal over $k$ iff $\tau_2\ne 0$. 
\item[(ii)] If $\tau_2\ne 0$, then for any $q$
\[
\dim_k D =\left\{
          \begin{array}{ll}
          2\lfloor q/3\rfloor & \text{for} \quad q\not\equiv 0\bmod 3\;,  \\
          {2\over 3}q - 1     & \text{for} \quad q\equiv 0\bmod 3\;.
          \end{array}
          \right.
\]
\end{description}
\end{lemma}
{\noindent \it Proof:\quad } 
(i)\quad If $q$ is a power of the characteristic of $k$, 
then $\tau_2 = (u-v)^{q-1}$, whence $\tau_2 \ne 0$ iff $u\ne v$ iff $C$
is nodal.
\smallskip\noindent
(ii) just evaluates the formula for $\dim_k D$ found above in terms of
$q\bmod 3$.$\hfill\square$
\[ \]

Putting everything together yields the Hilbert-Kunz 
function in the nodal case.

\begin{theorem} 
Let $C$ be a nodal cubic over a field $k$ of prime characteristic $p$.
For a power $q$ of $p$, the Hilbert-Kunz function at $q$ is 
\[ HK_C(q) =\left\{
            \begin{array}{ll} 
     {7\over 3}q^2 - {1\over 3}q - 1  & 
     \text{for} \quad q\not\equiv 2\bmod 3\; ,\\ 
     {7\over 3}q^2 - {1\over 3}q - {5\over 3} &
     \text{for} \quad q\equiv 2\bmod 3\;.
         \end{array}
         \right.
\]
\label{th:nodal}
\end{theorem}

If $C$ is a cuspidal cubic, then $\tau_2=0$ for any $q$ and 
we get immediately the generalized Hilbert-Kunz function --- 
in accordance with [8] and [1]:
\[ HK_{C,(x,y,z)}(q) =\left\{ 
                      \begin{array}{ll} 
             {7\over 3}q^2 & \text{for} \quad q\equiv 0\bmod 3\;,  \\
         {7\over 3}q^2 - {4\over 3} & \text{for} \quad q\not\equiv 0\bmod 3\;.
                      \end{array}
                      \right.
\]

Note however that,
if $q$ is not a power of the characteristic,
this last result will in general depend upon the
choice of the coordinate system made relative to the given
Weierstra\ss\ form.
The case of the generalized Hilbert-Kunz function for a nodal cubic
can be extracted as well --- 
and the dependence upon the coordinate system becomes 
apparent: If the distinct roots $u,v$ satisfy $u^q-v^q = 0$ for some $q$,
the generalized Hilbert-Kunz function ``jumps up'', 
it takes on the value from the cuspidal case. 
For any given $q$, we can avoid this situation by
replacing $y$ with $y+\alpha x$, for a general $\alpha\in k$. 
The curve $C$ is then still in Weierstra\ss\ form, and
with respect to $(x,y+\alpha x, z)$, 
the generalized Hilbert-Kunz function takes on the value predicted
by Theorem \ref{th:nodal}. 
Unless the algebraically closed field $k$ 
is an algebraic closure of a finite field, one 
can even find an $\alpha\in k$ that works for all $q$
simultaneously.

\section{Elliptic Curves in Odd Characteristic}
In this section, we prove the announced result for elliptic 
curves in odd characteristic and deduce that the Hilbert-Kunz 
multiplicity of a generic plane curve equals $\frac{3}{4}d$ when
$d\geq 2$.

\begin{theorem}  
Let $f(x,y,z)\in S=k[x,y,z]$ be a cubic polynomial defining
a plane elliptic curve over a field $k$ of odd characteristic
$p$. For any $n\in {\Bbb{N}}$ and $q=p^n$, 
the socle degree $a(q)$ of $\theta = S/(f+m^{[q]})$ is minimal,
\[a(q)=\frac {3}{2}q-\frac {1}{2}\quad,\] 
and the Hilbert-Kunz function of $R=S/(f)$ at $q$ is given by
\[ HK_R (q)=\frac {9}{4}q^2 - \frac {5}{4}\quad .\]
\label{theorem:main}
\end{theorem}

\newtheorem{corollary}{Corollary}
\begin{corollary}
For any field $k$ of prime characteristic $p$ 
and any integer $d\geq 2$,  
there is a curve $C\subset {\Bbb{P}}^2_k$
of degree $d$ whose Hilbert-Kunz multiplicity
achieves the minimum $\frac{3}{4}d\;$.
\end{corollary}

{\noindent \it Proof:\quad } As shown in [8], 
the Hilbert-Kunz multiplicity of the quadric
$g=x^2-yz$ equals $\frac{3}{2}$. For elliptic
curves in any prime characteristic, 
Theorem \ref{theorem:major} shows
that their Hilbert-Kunz multiplicities are minimal, 
equal to $\frac{9}{4}$.
As any integer $d\geq 2$ can be written 
$d=2u+3v$ for some $u,v\in {\Bbb{N}}$, 
additivity of the Hilbert-Kunz multiplicity, see [6],
implies that the curve of degree $d$, 
defined by $h=g^uf^v$, $f$ 
a nonsingular cubic, will achieve the minimum.
$\hfill \square$
\[ \]

{\noindent \it Remark 3. } 
Semi-continuity of
the Hilbert-Kunz multiplicity yields that the 
Hilbert-Kunz multiplicity of a
generic plane curve of degree $d\geq 2$ equals 
$\frac{3}{4}d$. Clearly
$c=1$ if the degree $d=1$. So the Hilbert-Kunz
multiplicity of a generic curve is rational and 
independent of the (positive) characteristic. 
\[ \]

In section 4.1, we recall
a classical result about determinants of Hankel
matrices whose entries are Legendre polynomials, and
in section 4.2,
we use it to determine the invariant
$a(q)$ and to establish Theorem \ref{theorem:main}.

\subsection {Hankel Determinants of Legendre Polynomials}             

The Hankel matrices associated
to a sequence $a=\{ a_i \}$  are 
\[ H^{(n)}_k (a) = \left(
                   \begin{array}{llll}
                   a_n & a_{n+1} & \cdots & a_{n+k-1}  \\
                   a_{n+1} & a_{n+2} & \cdots & a_{n+k} \\
                   \cdots & \cdots & \cdots & \cdots    \\
                   a_{n+k-1} & \cdots & \cdots & a_{n+2(k-1)}
                   \end{array}
                   \right)\quad,
\]
with corresponding Hankel determinants
\[  D^{(n)}_k (a) = \det H^{(n)}_k (a) \quad .  \]
The generating function for the Legendre polynomials,
$\{P_n (t)\}_{n\in\Bbb N},$ is  
\[  F(t,x)=\frac{1}{\sqrt{1-2tx+x^2}} =
\sum P_n (t) x^n \in{\Bbb Z}[\tfrac{1}{2},t][\![ x]\!] \quad .\]
For each $k$, consider the determinant of the following 
Hankel matrix whose entries are Legendre polynomials,
\[ D^{(0)}_k (P_i (t)) =\det \left(
              \begin{array}{llll}
              P_0 (t) & P_1 (t) & \cdots & P_{k-1} (t)  \\
              P_1 (t) & P_2 (t) & \cdots & P_k (t)     \\
              \cdots  & \cdots  & \cdots & \cdots      \\
              P_{k-1} (t) & \cdots & \cdots & P_{2k-2} (t)   
              \end{array}
              \right)\quad.
\]
In [2], J. Geronimus gave the following beautiful formula.

\begin{theorem}
\[ D^{(0)}_k (P_i (t)) = 
2^{- (k-1)^2} (t^2-1)^{\frac{1}{2} (k-1)k} \quad .  \]
\label{theorem:G}
\end{theorem}

P. Monsky communicated a direct proof to us that we now present.

\begin{lemma}
If $t\in {\Bbb{R}}$ and $t>1$, then
\[ P_n (t)=\frac{1}{\pi} \int^{\pi}_0 [t+\sqrt{t^2-1}
\cos\alpha]^n d \alpha   . \]
\end{lemma}

{\noindent \it Proof:\quad} 
If $x$ is small,
\[ (1-2tx+x^2)^{-\frac{1}{2}}
=[(1-tx)^2-(x\sqrt{t^2-1})^2]^{-\frac{1}{2}}     \]
\[=\frac{1}{\pi} \int^{\pi}_0 \frac{d\alpha}{(1-tx)-x\sqrt{t^2-1}
\cos\alpha }
=\frac{1}{\pi} \int^{\pi}_0 \frac{d\alpha}
{1-x[t+\sqrt{t^2-1}\cos\alpha ]}\quad . \]
The expansion of the integrand into a power series of $x$ 
yields the lemma. $\hfill \square$
\[ \]

{\noindent \it Proof of Theorem \ref{theorem:G} (P. Monsky):\quad}
As both sides of Geronimus' formula are polynomials in $t$,
it suffices to prove it for $t\in {\Bbb{R}}$ and $t>1$.

Let $V$ be the vector space of real
valued continuous functions on $[0,\pi]$. Define a symmetric
bilinear form $(\cdot,\cdot)$ on $V$ through
\[ (f,g)=\frac{1}{\pi} \int^{\pi}_0 f(\alpha)g(\alpha) d\alpha 
\quad . \]
If $h_1,\ldots,h_s \in V$, let
$\Delta (h_1,\ldots,h_s)$ be the determinant of the matrix
$((h_i,h_j))$. 
Set $g_m=[t+\sqrt{t^2-1} \cos\alpha]^m$. By
the preceding lemma, $(g_i,g_j) = P_{i+j}$. 
So the required Hankel determinant
is $\Delta (g_0,\ldots,g_{k-1})$. 

Let $V_m\subset V$ be the subspace spanned by
$\{1,\cos\alpha,\ldots,(\cos\alpha)^m\}$.
Then $f_m=(t^2-1)^{\frac{m}{2}} \cos (m\alpha) \in V_m$, 
and it is a linear combination of $g_0,\ldots,g_m$. 
Furthermore,
modulo $V_{m-1}$, $\cos(m\alpha) \equiv 2^{m-1} (\cos\alpha)^m$,
and consequently $f_m \equiv 2^{m-1} (t^2-1)^{\frac{m}{2}}
(\cos \alpha)^m \equiv 2^{m-1} g_m$. We conclude that
\[ \Delta (f_0,\cdots,f_{k-1})=(\prod^{k-1}_{m=1} 2^{m-1})^2
 \Delta (g_0,\cdots,g_{k-1})=2^{(k-2)(k-1)}\Delta (g_0,\cdots,g_{k-1}).\]
But using the orthogonality of the $f_i$ one finds that
\[ \Delta (f_0,\cdots,f_{k-1})=2^{-(k-1)}(t^2-1)^{\frac{(k-1)k}{2}} \]
and Theorem \ref{theorem:G} follows. $\hfill \square$
\[ \]

Now consider
\[ G(t,x)=\sqrt{1-2tx+x^2}=\sum \widetilde{P_n (t)} x^n\quad.  \]
As $F(t,x) G(t,x) =1$,
Geronimus' formula yields also the following corollary.

\begin{corollary}
\[ D^{(2)}_k (\widetilde{P_i (t)} )= (-1)^k   D^{(0)}_{k+1} (P_i (t)) 
             = (-2)^{-k^2} (t^2-1)^{\frac{1}{2} k(k+1)}
             \quad . \]
\end{corollary}

{\noindent \it Remark 4. }
The coefficients of Legendre polynomials are rational numbers
whose denominators are powers of $2$. Thus Geronimus' identity
and the above corollary hold over any ring in which $2$ 
is a unit, in particular over a field of odd characteristic.

\subsection{The Invariant $a(q)$ }

We first prove the following theorem showing that 
in odd characteristic there are no
nontrivial syzygies of low degree between 
the equation of an elliptic curve and Frobenius powers of the
variables.

\begin{theorem}
Let $k$ be a field of odd characteristic $p$,
and let $f\in k[x,y]$ be a cubic polynomial defining 
an elliptic curve in ${\Bbb{A}}^2_k$.
For any $q=p^n$, with $n\in {\Bbb{N}}$,
if $f|u x^q+v y^q+w$ for $u,v,w\in k[x,y]$ of degree at most
$\frac{1}{2}(q-1)$, then $f$ divides each of $u,v,w$.
\label{theorem:second}
\end{theorem}

{\noindent \it Proof:\quad} 
We give the proof for $q\equiv 1\bmod 4$.   
The argument in the other case, $q\equiv 3\bmod 4$, 
is analogous and left to the reader. Without loss
of generality assume that $k$ is algebraically closed.
Since the result
is invariant under the action of $GL(2, k)$, and
the characteristic of $k$ is odd,
we can put the cubic into the form $f=y^2-x(1-2tx+x^2)$  
with $t^2\neq 1$.

If $f | u x^q + v y^q + w $ for some $u , v , w $ 
of degree at most $\frac{1}{2} (q-1) $, then 
\begin{equation}
 u x^q + v y^q + w = f h \quad ,
\label{eq:main}
\end{equation}
where $h$ is a polynomial in $x,y$.
Set
\begin{eqnarray*}
 l & = &  \frac{q-1}{2}\quad,  \\ 
 g & = &  x(1-2tx+x^2)\quad .
\end{eqnarray*}
We can then write
\[ u=\sum^l_{i=0} a_i y^i 
    =A_0 + yA_1 + fu_1\quad,   \]
with $a_i\in k[x]$ of degree at most $l-i,
u_1 \in k[x,y]$ and
\[ A_0=\sum^{\frac{l}{2}}_{j=0} a_{2j}g^j\quad,\quad
   A_1=\sum^{\frac{l}{2}-1}_{j=0} a_{2j+1}g^j\quad, \] 
polynomials in $x$. Similarly, we write
\begin{eqnarray*}
   v & = & B_0 +yB_1+fv_1  \\
   w & = & C_0 +yC_1+fw_1  \\
 y^q & = & y^{2l+1}=yg^l+ f\gamma 
\end{eqnarray*}
for polynomials $v_1, w_1, \gamma\in k[x,y]; 
B_0,B_1,C_0,C_1\in k[x]$.
Equation (\ref{eq:main}) then becomes
\[ (x^qA_0+C_0+g^{l+1}B_1)+y(x^qA_1+C_1+g^lB_0)=fh_1 \quad .\]
Viewing both sides as polynomials in $y$, we get $h_1=0$ and
\begin{equation}
 x^qA_0+C_0+g^{l+1}B_1  =  0\quad ,
\label{eq:1}
\end{equation}
\begin{equation}  
 x^qA_1+C_1+g^lB_0   =  0\quad  .
\label{eq:2}
\end{equation}  
As $\deg C_0\leq \deg x^q A_0 \leq \frac{7}{2} l + 1$,
it follows that
\[ \deg B_1 \leq (\frac{7}{2} l + 1) - 
 3 (l + 1) = \frac{l}{2} -2\quad , \]
and similarly
$\deg B_0 \leq \frac{l}{2} -1$. 
Thus we can write 
\begin{eqnarray*} 
 B_1 & = & \alpha_{\frac{l}{2}-2} x^{\frac{l}{2}-2} + 
           \cdots + \alpha_0 \quad, \\ 
 B_0 & = & \beta_{\frac{l}{2}-1} x^{\frac{l}{2}-1}
           + \cdots + \beta_0\quad,    
\end{eqnarray*}  
for tuples
\[ \alpha = ( \alpha_i ) \in k^{\frac{l}{2}-1} \quad , 
\quad \beta= (\beta_i) \in k^{\frac{l}{2}} \quad .     \]
Since $\deg C_0 \leq \frac{3}{2} l $ and 
${\rm ord}(x^q A_0) \geq q=2l+1 $, 
the intermediate powers of $x$ in
$g^{l+1} B_1 $ have zero coefficients, whence we get 
a linear system of equations for $\alpha$, say
$E \alpha =0$, where
$E : k^{\frac{l}{2}-1} \longrightarrow  k^{\frac{l}{2}}$ 
is represented by the matrix
\[ E=\left(
     \begin{array}{llll}
     e_2 & e_3 & \cdots & e_{\frac{l}{2}}  \\
     e_3 & e_4 & \cdots & e_{\frac{l}{2}+1} \\
     e_4 & \cdots & \cdots & e_{\frac{l+2}{2}}  \\
     \cdots & \cdots & \cdots & \cdots    \\
     e_{\frac{l}{2}+1} & \cdots & \cdots & e_{l-1} 
     \end{array}
     \right)\quad,
\] 
whose entries $e_i$ are the coefficients in the expansion
\[ (1-2tx+x^2)^{\frac{q+1}{2}} = x^{q+1} + e_q x^q +\cdots
  + e_0 \quad . \]
Analogously, the corresponding powers of $x$ in $g^l B_0$ 
yield a system of equations for $\beta$, say
$H \beta = 0$,
where
$H : k^{\frac{l}{2}} \longrightarrow k^{\frac{l}{2}+2}$ is
represented by the matrix
\[ H=\left(
    \begin{array}{llll}
    h_0 & h_1 & \cdots & h_{\frac{l}{2}-1}  \\
    h_1 & h_2 & \cdots & h_{\frac{l}{2}}  \\
    \cdots & \cdots & \cdots & \cdots    \\
    h_{\frac{l}{2}+1} & \cdots & \cdots & h_l  
    \end{array}
    \right)\quad,  
\]
whose entries $h_i$ are the coefficients in the expansion 
\[ (1-2tx+x^2)^{\frac{q-1}{2}}= x^{q-1} + h_{q-2} x^{q-2} +
  \cdots + h_0\quad . \]    
As $q$ is a power of the characteristic $p$, one has
\[(1-2tx+x^2)^{\frac{q-1}{2}} \equiv 
\frac{(1 - (2tx)^q + x^{2q})^{\frac{1}{2}}}
{\sqrt{1-2tx+x^2}} \equiv (1 - (2tx)^q + x^{2q})^{\frac{1}{2}}
\sum P_n (t) x^n \bmod p\quad,\]
whence $h_i \equiv P_i (t)\bmod p$ and, analogously,
$e_i \equiv \widetilde{P_i (t)}\bmod p$ for $i<q$, where
$P_i (t)$ and $\widetilde{P_i (t)}$ are as in section 4.1. 
As $t^2\neq 1$ by assumption, Geronimus' Theorem and its Corollary
imply
\[ {\rm rank} \; E = \frac{l}{2} -1 \quad,\quad
 {\rm rank} \; H = \frac{l}{2}\quad , \]
whence each $\alpha_i$ or $\beta_j$ equals zero, thus
$B_0=B_1=0$, so that $v=fv_1$.
As $\deg C_i <{\rm ord} (x^qA_i)$; for $i=0,1$;
it follows further from equations (\ref{eq:1}) and (\ref{eq:2})
that $C_i=A_i=0$ and the theorem follows. 
$\hfill \square$
\[ \]

Now we can finish the {\it Proof of Theorem \ref{theorem:main}:}

As $d=3<3(q-1)$, for any power $q=p^n; n\in {\Bbb{N}};$ 
of an odd prime $p$,
Theorem \ref{theorem:minimal} yields the lower bound
$a(q) \geq \frac{3}{2}q-\frac{1}{2}$, 
and the upper bound $\imath (q)\leq \frac{3}{2} q-\frac{5}{2}$.
It remains thus to show
$\vartheta_{\frac{3}{2} q-\frac{7}{2} }=0$,  
or, equivalently, 
if  $f | u x^q+v y^q + w z^q $ for 
$u ,v ,w \in k[x,y,z]_{\frac{1}{2} (q-1) }$, then
$f | u ,v , w $.
As it suffices to verify the above statement in
the affine part $(z=1)$ of ${\Bbb{P}}^2_k $, the result in
Theorem \ref{theorem:second} finishes the proof. 
$\hfill \square$
          
\section{Cayley's Cubic Surface}

Let $S=k[x,y,z,w]$ be the polynomial ring in four variables over
an arbitrary field $k$ and let $f=xyz+xyw+xzw+yzw$ be the 
Cayley cubic.
We consider the generalized Hilbert-Kunz function of $R=S/f$,
given at $q\in\Bbb N$ through
\[  HK_{R,(x,y,z,w)} (q)=\dim_k S/(f,x^q,y^q,z^q,w^q) \quad. \]

\begin{theorem}
The socle degree of the artinian 
ring $\theta= S/(f,x^q,y^q,z^q,w^q)$ is
\[ a(q)=\left\{
     \begin{array}{ll}
     0    & {\text if} \quad q=1\;, \\
     2n-1 & {\text if} \quad q>1\;,
     \end{array}
     \right.
\]
and the value of the generalized Hilbert-Kunz function 
of Cayley's cubic at $q\in\Bbb N$ is 
\[  HK_{R,(x,y,z,w)} (q)= 2q^3 - q\quad.\]
\label{theorem:third}
\end{theorem}
{\noindent \it Proof:\quad } 
If $q=1$, then $\theta \cong k$ and $a(1)=0$.
Now assume $q>1$. 
Since $d=3<4(q-1)$, Theorem \ref{theorem:minimal}
yields the lower bound $a(q)\geq 2q-1$. 
Thus it remains to show $\theta_{2q}=0$, i.e,
that any monomial $x^iy^jz^kw^l\in \theta_{2q}$
is equivalent to $0$.
We prove this by ``descent'' on the 
{\it dominant exponent\/}, $e=\max\{i,j,k,l\}$,
of a monomial $x^iy^jz^kw^l$ in $S_{2q}$.

Any monomial with sufficiently large dominant exponent ---
for example, when it exceeds $q$ --- will be 
equivalent to $0$. We may thus assume
that for a fixed integer $e$, any monomial whose
dominant exponent exceeds $e$ is equivalent to $0$,
to consider then monomials with dominant exponent equal to $e$. 
Due to symmetry, it suffices to consider monomials $x^iy^jz^kw^l$ 
with $e=i\geq j \geq k \geq l$.

{\it Case 1: \quad} 
If $l>0$, then 
\begin{eqnarray*}
 x^iy^jz^kw^l & \equiv & 
 x^iy^{j-l}z^{k-l} \cdot (-1)^l x^l (yz+yw+zw)^l \bmod f\\ 
              & \equiv & 
 (-1)^lx^{i+l} y^{j-l}z^{k-l}(yz+yw+zw)^l \bmod f\\ 
              & \equiv & 0\bmod (f,x^q,y^q,z^q,w^q) 
\end{eqnarray*}
by assumption.  
 
{\it Case 2: \quad}
Now suppose $l=0$. 
If $k\leq 1$, then $i+j+1\geq 2q$ with $i\geq j$, so that
$i\geq q$ whence the monomial is obviously equivalent to
$0$. If $k>1$,
consider first the monomial $g=x^{i-k+1}y^jz^kw^{k-1}$
and set $h=\min\{ i-k+1,k-1 \}$. Since $j>i-k+1>0$ and  
$j\geq k>k-1>0$, the argument employed above shows that $g$ is
equivalent to a linear combination of monomials whose dominant
exponent is $j+h$.
As $j+h>i$, our assumption insures that 
$g$ is equivalent to $0$. On the other hand,
\[ g=x^{i-k+1}y^jz^kw^{k-1}
    \equiv  (-1)^{k-1} x^i y^{j-k+1} z (yz+yw+zw)^{k-1}
    	\bmod f\quad ,\]  
and in the expansion of the right hand side of this congruence,
each term involving $w$ involves all four variables,
has dominant exponent $e=i$, and is thus equivalent to
$0$ by Case 1.
So
\[  g=x^{i-k+1}y^jz^kw^{k-1}
    \equiv  (-1)^{k-1} x^i y^{j-k+1}z (yz)^{k-1}
    \equiv (-1)^{k-1} x^i y^j z^k\quad,\]
whence $x^iy^jz^k\equiv 0$.

The claimed result for the generalized Hilbert-Kunz function 
follows now from Theorem \ref{theorem:minimal}, as $d=3 < 4(q-1)$
whenever $q >1$, and its validity for $q=0,1$ is clear.
$\hfill \square$
\[ \]

\begin{corollary}
For any field $k$ of prime characteristic $p$ 
and any integer $d\geq 2$,  
there is a surface $X\subset {\Bbb{P}}^3_k$
of degree $d$ whose Hilbert-Kunz multiplicity
achieves $\frac{2}{3}d$, the minimum possible for such surfaces.
\end{corollary}

{\noindent \it Proof:\quad} 
The Hilbert-Kunz multiplicity
of the quadric surface $g=xy-zw$ equals $\frac{4}{3}$
by [1].
As just established, 
the Hilbert-Kunz multiplicity of the 
Cayley cubic $f$ is equal to $2$. Since
$d=2u+3v$ for some $u,v\in {\Bbb{N}}$, additivity of the 
Hilbert-Kunz
multiplicity implies that the surface defined by $g^uf^v$ 
has Hilbert-Kunz multiplicity equal to $\frac{2}{3}d$.   
$\hfill \square$
\[ \]

{\noindent \it Remark 5. } 
For any field $k$ of positive characteristic,
by virtue of the above corollary and semi-continuity, 
a generic surface in 
${\Bbb{P}}^3_k$ of degree $d\geq 2$ achieves the minimal 
Hilbert-Kunz multiplicity $\frac{2}{3}d$. 
Also, $c=1$ if $d=1$.
So the Hilbert-Kunz multiplicity of a generic surface 
is rational and independent of the characteristic. 
Note that the Hilbert-Kunz multiplicity of
Cayley's cubic is minimal although this surface is
singular --- in contrast to the case of cubic curves.
\[ \]

{\noindent \it Acknowledgments.\quad}  
The authors would like to thank P. Monsky, whose influence on
this paper should be clear to every reader. We are
especially grateful for his permission to include the
direct proof of Geronimus' theorem and for
modifications that lead to a shorter proof 
of Theorem \ref{theorem:second}. 
We also thank K.~Pardue, who got
us interested in the subject, 
and A.~Conca for useful conversations.

\end{document}